\documentclass[11pt]{amsart}

\usepackage{amssymb,url}
\usepackage[british]{babel}

\theoremstyle{remark}     
\newtheorem{rmk}{Remark}[section]\newtheorem*{rmk*}{Remark}
\newtheorem*{ex*}{Example}
\newtheorem*{acknowledgements}{Acknowledgements}

\swapnumbers              
\theoremstyle{plain}      
\newtheorem{lemma}[rmk]{Lemma}
\newtheorem{prop}[rmk]{Proposition}
\newtheorem{thm}[rmk]{Theorem}
\newtheorem{cor}[rmk]{Corollary}

\theoremstyle{definition} 
\newtheorem{defn}[rmk]{Definition}

\numberwithin{equation}{section}

\newcommand{\hodge}{{*}}
\newcommand{\lie}[1]{\mathfrak{#1}}
\newcommand{\Lie}[1]{\textsl{#1}}
\newcommand{\G}{G_2}
\newcommand{\SU}{\operatorname{\Lie{SU}}(3)}

\newcommand{\g}{\lie{g}}

\DeclareMathOperator{\Ker}{Ker}
\DeclareMathOperator{\im}{Im}
\newcommand{\W}{\mathcal{W}}
\newcommand{\X}{\mathcal{X}}

\newcommand{\psp}{\psi^+}
\newcommand{\psm}{\psi^-}
\newcommand{\pspm}{\psi^\pm}

\newcommand{\iso}{\cong}
\newcommand{\then}{\Longrightarrow}
\newcommand{\lto}{\longrightarrow}

\newcommand{\lan}{\langle}
\newcommand{\ran}{\rangle}
\newcommand{\hook}{\lrcorner\,}

\newcommand{\rcomp}[1]{\mathopen{\big[\mkern-5mu\big[} #1
\mathclose{\big]\mkern-5mu\big]}}
\newcommand{\rreal}[1]{\bigl[#1\bigr]}

\newcommand{\ol}{\overline}
\newcommand{\sym}{\mathcal{S}}
 
\newcommand{\oo}{\omega^2}
\newcommand{\1}{e^1}\newcommand{\2}{e^2}\newcommand{\3}{e^3}
\newcommand{\4}{e^4}\newcommand{\5}{e^5}\newcommand{\6}{e^6}
\renewcommand{\geq}{\geqslant}

\begin{document}
\title[$G_2$T-structures from nilmanifolds]{$\boldsymbol{G_2}$-structures
with torsion from\\ half-integrable nilmanifolds} 
\author{Simon G.~Chiossi}
\address[Chiossi]{Department of Mathematics and Computer Science,
University of Southern Denmark, Campusvej 55, 5230 Odense M, Denmark}
\email{chiossi@imada.sdu.dk}
\author{Andrew Swann}
\address[Swann]{Department of Mathematics and Computer Science,
University of Southern Denmark, Campusvej 55, 5230 Odense M, Denmark}
\email{swann@imada.sdu.dk}

\begin{abstract}
  The equations for a $G_2$-structure with torsion on a product $M^7 = N^6
  \times S^1$ are studied in relation to the induced
  $\Lie{SU}(3)$-structure on~$N^6$.  All solutions are found in the case
  when the Lee-form of the $G_2$-structure is non-zero and $N^6$ is a
  six-dimensional nilmanifold with half-integrable $\Lie{SU}(3)$-structure.
  Special properties of the torsion of these solutions are discussed.
\end{abstract}


\maketitle

\section{Introduction}
Connections with torsion have been objects of geometrical study for many
years. Interest in this subject has been increased by considerations from
supersymmetric string- and M-theories~\cite{Strominger:superstrings}, where
connections with skew-symmetric torsion coming from $G$-structures
distinguished by spinors play an important r\^ole.  Recent mathematical
discussion of such connections may be found, for example, in
\cite{Friedrich-I:skew,Gibbons-PS:hkt-okt,Grantcharov-P:HKT}, where
Hermitian manifolds, $\G$, $\operatorname{\Lie{Spin}}(7)$ and quaternionic
geometries occur.  Particular importance is attached to geometries in
dimensions $6$ and~$7$ given by $\SU$- and $\G$-structures.

In this paper we investigate the explicit construction of $\G$ geometries
with torsion ($\G$T-structures) from products with six-dimensional
$\SU$-manifolds.  The initial data on $N^6$ is an almost Hermitian
structure $(J,h,\omega)$ together with a distinguished complex volume
$\Psi=\psi^++i\psi^-$.  A $\G$-structure may then be built on the product
$M^7=N^6\times S^1$ by using the three-form $\varphi = \omega\wedge
dt+\psi^+$.  

Following Gray \& Hervella \cite{Gray-H:16}, one approach to the study of
metric $G$-structures in general is via consideration of the components of
the intrinsic torsion in the irreducible $G$-modules of
$T^*M\otimes\mathfrak g^\bot$.  For~$\G$, this space splits into four
components $\X_i$, $i=1,\dots,4$ \cite{Fernandez-G:G2} and the conditions
on $\varphi$ for the intrinsic torsion~$\tau^{(M)}$ to lie in a given
combination of these spaces have been determined explicitly
in~\cite{Cabrera:G2}.  The pair $(M^7,\varphi)$ is a $\G$T-structure
exactly when the intrinsic torsion has no component in $\X_2\iso\mathfrak
g_2$, i.e.,
\begin{equation}
  \label{integrable}
  \tau^{(M)} \in \X_1\oplus \X_3\oplus\X_4.
\end{equation}
Slightly misleadingly, such structures are sometimes referred to as
`integrable $\G$-structures' in the literature, despite the fact that the
Riemannian holonomy need not reduce.  Special cases of this geometry are
studied in \cite{Cleyton-I:cG2, Bryant:G2, Friedrich-KMS:G2,
Fino:conformal-G2}.

The corresponding refinement of the Gray-Hervella classification for
$\SU$-\hspace{0pt}structures was computed in~\cite{Chiossi:PhD}.  Our first
task, in~\S\ref{sec:structures}, is to relate the two decompositions on
$N^6$ and $M^7$ concentrating on the situation for $\G$T-structures and
refining results of~\cite{Chiossi-S:SU3-G2}.

We then turn to consideration of particular examples.  In
\cite{Fino-PS:SKT}, six-dimensional nilmanifolds were successfully used to
give examples of SKT structures: KT geometry is `K\"ahler with torsion' and
consists of a Hermitian manifold $(N,J,h,\omega)$ together with its Bismut
connection, the unique Hermitian connection with totally skew-symmetric
torsion (essentially $Jd\omega$), see~\cite{Gauduchon:Hermitian-Dirac};
the structure is `strong' (SKT) when the torsion is a closed form, i.e.,
$dJd\omega=0$.  A nilmanifold is a compact quotient of a nilpotent Lie
group.  In dimension~$6$, each nilpotent Lie algebra has a basis with
rational structure coefficients and so~\cite{Malcev:rational} the
corresponding $1$-connected group~$G$ admits a co-compact lattice~$\Gamma$.
Much is known about the topology \cite{Nomizu:cohomology} and geometry of
these manifolds.  In particular, they do not admit K\"ahler metrics unless
$G$~is Abelian~\cite{Benson-G:nil, Cordero-FG:symplectic-not-Kaehler,
Campana:nil-Kaehler}, a fact related to the theory of minimal
models~\cite{Deligne-GMS:homotopy, Hasegawa:minimal}, and several authors
have studied their complex, Hermitian and symplectic geometry
\cite{Abbena-GS:aH-nil, Console-F:Dolbeault, Snow:complex-reductive,
Salamon:complex-nil, Cordero-FU:pK-6}.

We thus study the case of $N^6=\Gamma\backslash G$ with an $\SU$-structure
that pulls back to an invariant $\SU$-structure on~$G$.  However, instead
of looking at the full system of equations for a $\G$T-structure on
$M^7=N^6\times S^1$, we consider a special case when $N^6$ is
\emph{half-integrable}, see~\S\ref{sec:integrable}.  Restricting to the
case where $J$~is integrable is too severe, we merely obtain structures
with $\tau^{(M)} \in \X_1\oplus\X_4$.  The half-integrability condition
considered here and in~\cite{Chiossi-S:SU3-G2} is a weaker restriction on
the $\SU$-structure that is of interest in its own right: these are exactly
the $\SU$-structures that appear as hypersurfaces in manifolds of
holonomy~$\G$ (Joyce manifolds) and indeed Hitchin~\cite{Hitchin:3-form,
Hitchin:forms} has shown how the holonomy metric may be obtained by a flow
on the space of $\SU$-structures.

In \S\ref{HLNLAS}, we give a full classification of the invariant
half-integrable nilmanifolds $N^6 = \Gamma\backslash G$ and their
$\SU$-structures such that $M^7=N^6\times S^1$ carries a $\G$T-structure.
The proof relies on a detailed study of compatibility between the nilpotent
algebra structure and the $\SU$-geometry and is facilitated by
consideration of a complex version of the equations.  A key ingredient to
finding concrete solutions is the classification of six-dimensional
nilpotent Lie algebras as presented in~\cite{Salamon:complex-nil}.  In our
situation we find six different non-trivial Lie algebras from the
classification; in each case the relevant half-integrable $\SU$-structures
on $G$ are parameterised by at most two essential variables,
see~Theorem~\ref{output}.  In fact, these six algebras are closely related
to each other, and all are obtainable as degenerations of a single example.

We then analyse in \S\ref{SGT} the special properties of the derivative
of~$T$ for the $\G$-structures obtained.  We note examples where the
K\"ahler form $\omega$ on~$N^6$ is a eigenform for the Laplacian. The
eigenform property of $\omega$ is shared by `strong' $\G$T-structures of
the type suggested by the physical literature, but for a different
eigenvalue.  Interestingly, one of our examples also occurs
in~\cite{Ivanov-I:instantons} as an example of an `instanton'.  Future work
will concentrate on determining properties of the $\G$-holonomy metrics
containing these $\SU$-structures on hypersurfaces.

\begin{acknowledgements}
  S.G.C.~thanks Simon Salamon for the interest shown, and the
  Department of Mathematics and Computer Science in Odense where this paper
  was prepared. His gratitude also goes to Thomas Friedrich and Ilka
  Agricola for the invitation to Humboldt Universit\"at, where partial
  results were presented.
  
  Both authors are members of the \textsc{Edge} Research Training Network
  \textsc{hprn-ct-\oldstylenums{2000}-\oldstylenums{00101}}, supported by
  the European Human Potential Programme.  They thank Richard Cleyton for
  useful discussions in the early stages of this project.
\end{acknowledgements}

\section{$\G$-structures with torsion}
\label{sec:structures}
\noindent Denote by $N$ a manifold of (real) dimension 6 with a
$\operatorname{\Lie U}(3)$-structure and trivial canonical bundle.  $N$ is
thus equipped with an orthogonal almost complex structure $J$ and a
non-degenerate $2$-form $\omega$. The induced Riemannian metric $h$
distinguishes the circle consisting of elements of unit norm in the
$2$-dimensional space $[\![\Lambda^{3,0}]\!]$, and an $\SU$-structure is
determined by the choice of a real $3$-form $\psp$ lying in this
$S^1$-bundle at each point. The associated $(3,0)$-form is
$\Psi=\psp+i\psm$ with $\psm=J\psp$, and the knowledge of the tensors
$J,\omega,\Psi$ determines the geometry of the manifold in full, though
only $\omega,\psp$ are, strictly speaking, necessary to pin down the
reduction to $\SU$.

Due to the local nature of the set-up, all descriptions will be valid at
least pointwise, so by choosing an orthonormal basis $e^1,\dots , e^6$ for
the cotangent bundle $T^*N$, we define the following forms
\begin{equation}\label{stdsus}
  \begin{gathered}
    \omega=e^{12}+e^{34}+e^{56}\in\Lambda^{1,1}N,\\
    \psp+i\psm=(e^1+ie^2)\wedge(e^3+ie^4)\wedge(e^5+ie^6)\in\Lambda^{3,0}N,
  \end{gathered}
\end{equation} 
relative to the decomposition into types determined by $J$. We will freely
use the familiar notation $e^{ij\dots}$ to indicate $e^i\wedge
e^j\wedge\dots$, so for instance $\psp=e^{135}-e^{146}-e^{236}-e^{245}$ and
$\psm=e^{136}+e^{145}+e^{235}-e^{246}$. In an even more concise way we
shall sometimes write, say the K\"ahler form, as $\omega=12+34+56$.

The classical Gray-Hervella decomposition of the Hermitian intrinsic
torsion space into four irreducible representations $\W_k,\ k=1,\dots , 4$,
was extended to tackle $\SU$ reductions in~\cite{Chiossi-S:SU3-G2} (but see
also~\cite{Bor-HL:Bochner}). In the former some relations between $\G$
manifolds and underlying $\SU$-structures were taken into consideration,
and we shall adhere to the same notation throughout. While retaining curly
symbols for the almost Hermitian classes, $W_j$ denotes the corresponding
intrinsic torsion component. The new elements in the theory are the
presence of an extra fifth class $\W_5$ whose component depends on
$(d\psp)^{3,1}$, and the reducibility of the $U(3)$-modules $\W_1,\W_2$
which split into symmetric halves, denoted by $\W_1^\pm,\W_2^\pm$. This
allows one to introduce a correspondence between hypersurfaces of
$\G$-manifolds and $S^1$-quotients, a fact reflected in the self-duality
versus anti-self-duality picture in 4 dimensions~\cite{Chiossi-S:SU3-G2}.

We consider a Riemannian product $M^7$ of $N$ with a circle, endowed with
product metric $g$. By means of the almost Hermitian structure the manifold
$M=N\times S^1$ naturally inherits a $\G$-structure by declaring
\begin{equation}\label{phi}
  \varphi=\omega\wedge dt +\psp\in\Lambda^3T^*M.
\end{equation}
In the terminology of~\cite{Hitchin:forms} this three-form is stable and
defines a reduction of the structure bundle to the exceptional group $\G$.
Though relying on standard
references~\cite{Bryant:exceptional,Salamon:holonomy,Joyce:book} for the
theory of $\G$-structures, we recall here the fundamental fact that a form
of the kind \eqref{phi} completely specifies the Riemannian metric $g$, an
orientation and Hodge operator $*$. The basis $e^1,\dots , e^6,e^7=dt$ is
also orthonormal for the metric $g$ determined by the inclusion $\G\subset
SO(7)$, and one has $\varphi=127+347+567+135-146-236-245$.

It is an easy matter of calculation to see that condition \eqref{integrable}
holds if and only if
\begin{equation}
  \label{g2teqn}
  d\hodge\varphi=\theta\wedge\hodge\varphi
\end{equation}
for some $1$-form $\theta$ (see~\cite{Cabrera:G2}), poignantly called by
Friedrich and Ivanov the Lee form of the $\G$-structure~\cite{Friedrich-I:Ks7}. 
We shall be mainly concerned with the case
where $\theta$ does not vanish. The particular interest of this class is
clear by the result that 

\begin{thm}
  {\rm\cite{Friedrich-I:skew}} On a $\G$-manifold $(M,\varphi)$ the
  following are equivalent:
  \begin{enumerate}
  \item[{\sl 1)}] the torsion of the $\G$-structure is a three-form
  $T\in\Lambda^3T^*M$; 
  \item[{\sl 2)}] there exists a unique linear connection with skew-symmetric
    torsion
    \begin{equation}
      \label{torsion}
      T=\tfrac{1}{6}\lan d\varphi,\hodge\varphi\ran\varphi - \hodge
      d\varphi+\hodge(\theta\wedge\varphi),
    \end{equation}
    where $\lan\ ,\ \ran$ is the inner product given by the metric.
  \end{enumerate}
\end{thm}

Reflecting the orthogonal splitting $T^*M=\mathbb R^6\oplus\mathbb R$, we
decompose
\begin{equation}
  \label{Lee}
  \theta=\beta+\lambda\,dt
\end{equation}
for some $\lambda\in\mathbb R$ and $1$-form $\beta\in\Lambda^1N^6$.  One
immediate consequence of \eqref{g2teqn} is that the Lee form is closed, for
the wedging map with $\hodge\varphi$ is a one-to-one homomorphism between
$\Lambda^2=\mathbb R^7\oplus\g_2$ and $\Lambda^6=\Lambda^1=\mathbb R^7$.
The effects of this fact in physics are known, and will be recalled later.
The class of $\G$-manifolds with closed Lee
form is known as that of locally conformally balanced $\G$-structures.
The term `balanced' reflects the six-dimensional set-up with the same
name, where balanced, or cosymplectic, refers to a Hermitian structure with
$\vartheta=-Jd^*\omega=0$ (which is called Lee form, too). In terms of
intrinsic torsion $\tau$, this amounts to $\tau^{(N)}\in\W_3$, echoed
in seven dimensions by co-calibrated $\G$-structures, for which
$\tau^{(M)}$ belongs to $\X_1\oplus\X_3$.

We begin investigating the geometric properties of \eqref{g2teqn} by
expanding the exterior derivative of $\hodge\varphi=\psm\wedge
dt+\tfrac12\omega\wedge\omega$ into
\begin{equation*}
  d\hodge\varphi=d\psm dt+\omega d\omega
  =\beta\psm dt+\tfrac12\beta\oo+\tfrac12\lambda\oo dt
  =(\lambda dt+\beta)\hodge\varphi,
\end{equation*} 
where we start dropping wedge product signs to lighten expressions.
Comparing components yields
\begin{equation}
  \label{g2t}
  \left\{
    \begin{aligned}
      d\psm &=\beta \psm+\tfrac12\lambda\,\oo,\\
      \omega d\omega&=\tfrac12\beta\oo, 
    \end{aligned}
  \right.
\end{equation}
which we shall refer to as the \textbf{$\boldsymbol{G_2T}$ equations} for
$N$.  Since these are equivalent to $\X_2=0$, we obtain the first
restrictions

\begin{lemma}
  Whenever the $\G$-structure of $M=N\times S^1$ has a three-form torsion,
  the intrinsic torsion of $N$ satisfies
  \begin{equation}
    \label{g2tclasses}
    W_2^-=0,\qquad W_1^-=\tfrac12\lambda,\qquad W_5=-2W_4.
  \end{equation}
\end{lemma}

\begin{proof}
  an immediate consequence of~\cite{Chiossi-S:SU3-G2}.
\end{proof}

\begin{rmk*}
  The vanishing of the linear combination $2W_4+W_5$ is a requirement in
  four-dimensional ${\mathcal N}=1$ spacetime supersymmetry to give rise to
  supersymmetric compactifications of heterotic string
  theory~\cite{Strominger:superstrings}.
\end{rmk*}

The $\SU$ components $W_4,W_5$ are explicitly given by:
\begin{equation*}
  W_4=\tfrac12\omega\hook d\omega=\tfrac{1}{4}\,\beta\quad \text{so that}\quad
  W_5=-\tfrac12\beta.
\end{equation*}
Let us decompose $d\psp$ into types and set
$d\psp=\gamma\psp+W_2^+\omega+W_1^+\oo$. The definition of
$W_5=\tfrac12\psp\hook d\psp$ tells us that $\gamma=-\beta$, so
\begin{equation*} 
  d\psp=-\beta\psp+W_2^+\omega+W_1^+\oo.
\end{equation*}
We now repeat the procedure for the derivative of the K\"ahler form,
initially prescribing
$d\omega=\tfrac12\beta\omega+\Omega+\nu^+\psp+\nu^-\psm$, for
$\nu^\pm\in\mathbb R$. Here $\Omega=(d\omega)^{2,1}_0$ represents the
component in $\W_3$. In local coordinates, both $\omega^3$ and
$\psp\wedge\psm$ are multiples of the volume form:
\begin{equation*}
  \omega^3=6e^{12\dots 6},\qquad \psp\psm=4e^{12\dots 6}
\end{equation*}
in agreement with the compatibility equation
$\psp\psm=\tfrac{2}{3}\omega^3$ (cf.~\cite{Hitchin:3-form}).  Using the
definition of $W_1^\pm$, we find
$W_1^\pm\,\omega^3=d\pspm\wedge\omega=\pspm\wedge
d\omega=6W_1^\pm\,e^{12\dots 6}$, and $\pspm\wedge\omega=0$ follows by type
considerations, whence
\begin{alignat*}{3}
  \psp\wedge d\omega&=\nu^-\psp\psm&\quad&\then\quad&\nu^-
  &=\tfrac{3}{2}\,W_1^+\\
  \psm\wedge d\omega&=\nu^+\psm\psp&&\then&\nu^+&=-\tfrac{3}{4}\,\lambda.
\end{alignat*}
At last then, we are able to express
\begin{equation*}
  d\omega=\tfrac12\beta\omega+\omega-\tfrac{3}{4}\lambda\psp
  +\tfrac{3}{2}W_1^+\psm.
\end{equation*}

Aiming at the expression for the totally skew-symmetric torsion, we
determine each term of \eqref{torsion} separately, starting by computing
\begin{equation*}
  \begin{split}
    d\varphi=d\omega dt+d\psp
    & =\tfrac12\beta\omega dt + \omega dt + \tfrac{3}{2}W_1^+\psm
    dt-\tfrac{3}{2}W_1^-\psp dt\\ 
    & \qquad -\beta\psp+W_2^+\omega+W_1^+\oo,
  \end{split}
\end{equation*}
whence $\lan d\varphi,\hodge\varphi\ran=\tfrac{3}{2}W_1^+\|\psm\|^2+\tfrac12
W_1^-\|\psp\|^2=12W_1^+$.

\begin{rmk}\label{KS}
  This last term has a relevant physical meaning: its vanishing is
  precisely the Killing spinor equation $\nabla\eta=0$ prescribing the
  existence of a parallel spinor field $\eta$ with respect to the torsion
  connection (consult~\cite{Friedrich-I:Ks7,Strominger:superstrings}).
  There is a second equation akin to this, namely $\theta=-2d\phi$, where
  $\phi$ represents the dilation function of string theory. We have already
  shown that $\theta$ is closed (so locally exact) in the present setting,
  see \eqref{Lee}.
\end{rmk}

Paying a little attention to the different behaviour of the Hodge star
operators, denoted $*_6$ if acting on $\mathbb R^6$, we have
\begin{equation*}
  \begin{split}
  *d\varphi&=\tfrac12{*}(\beta\omega dt)+*(\omega
   dt) + \tfrac{3}{2}W_1^-\psm + \tfrac{3}{2}W_1^+\psp-*(\beta\psp) +
   *(W_2^+\omega)+2W_1^+\omega dt\\
   &=-\tfrac12{J\beta}\wedge\omega+*_6\omega + \tfrac{3}{2}W_1^-\psm +
   \tfrac{3}{2}W_1^+\psp+\beta^\sharp\hook\psm 
   dt+*(W_2^+\omega)+2W_1^+\omega dt.
 \end{split}
\end{equation*}
Eventually, using the $\G$ Lee form \eqref{Lee},
\begin{equation*}
  \hodge(\theta\wedge\varphi)
  =\lambda\psm-J\beta\wedge\omega-\beta^\sharp\hook\psm dt.
\end{equation*}
Collecting the relevant terms we rewrite \eqref{torsion} as follows:
\begin{equation*}
  T=-\tfrac12 J\beta\wedge\omega -2\beta^\sharp\hook\psm
  dt-*_6\Omega-*(W_2^+\omega)+\tfrac12 
  W_1^-\psm+\tfrac12 W_1^+\psp,
\end{equation*}
and expressing $*_6\Omega$ in terms of $d\omega$ we arrive at the more
concise
\begin{equation}\label{T}
  T=-2\beta^\sharp\hook\psm
  dt-*(W_2^+\omega)-*_6d\omega+2W_1^+\psp-\tfrac12\lambda\psm.
\end{equation} 
To have more readable formul{\ae} we shall denote the $1$-form $\beta$ and
its dual vector $\beta^\sharp$ by the same symbol.

Furthermore, we shall refer to
\begin{equation}\label{*T}
  \hodge T=2\beta\psp-W_2^+\omega+d\omega dt+2W_1^+\psm
  dt+\tfrac12\lambda\psp dt 
\end{equation}
in the sequel, and to the differentials
\begin{equation}\label{dTd*T}
  \begin{split}
    dT
    &= -d{*}_6d\omega - 2W_1^+\beta\psp + 2W_1^+W_2^+\omega -
    \tfrac12\lambda\beta\psm+\bigl(2(W_1^+)^2-\tfrac{1}{4}\lambda^2\bigr)\oo\\  
    &\qquad-2d(\beta\hook\psm)dt-d{*}(W_2^+\omega)\\
    d\hodge T
    &=-2d(W_2^+\omega) + \tfrac12\lambda W_2^+\omega dt + 2W_1^+\beta\psm dt -
    \tfrac12\lambda\beta\psp dt+\tfrac{3}{2}\lambda W_1^+\oo dt\\ 
    &\qquad-2\beta W_2^+\omega-2W_1^+\beta\oo.
  \end{split}
\end{equation}

\begin{ex*} 
  Let $N$ be the Iwasawa manifold, that is, the compact quotient of the
  complexified Heisenberg group $G={\mathcal{H}_3}^\mathbb C$ by the sublattice of
  matrices with Gaussian integer entries. This nilmanifold arises from the
  Lie algebra $\g$ of $G$ given by
  \begin{equation*}
    de^i=(0,0,0,0,13+42,14+23).
  \end{equation*}
  Endow $N$ with the non-integrable almost complex structure $J_3$, in the
  notation of~\cite{Abbena-GS:aH-nil}, and select the following invariant
  $\SU$-structure
  \begin{equation*}
    \omega_3=-12-34+56,\qquad\psp=135+146+236-245.
  \end{equation*}
  Whilst the $\G$T equations reduce to $ d\psm=0, d\oo=0$, all geometric
  information is determined by
  \begin{equation*}
    d\psp=4e^{1234}\in\Lambda^{2,2},\qquad d\omega_3=\psm.
  \end{equation*}
  It is known that the $\G$-structure $\varphi=\omega_3\wedge e^7+\psp$ is
  co-calibrated, and in fact $d\hodge\varphi=\theta\wedge\hodge\varphi=0$
  implies $\theta=0$, thus $\beta=0,\lambda=0$.
  
  For later purposes, notice that the three-form torsion is given by
  $T=\tfrac{2}{3}\varphi-4e^{567}$, hence in particular is not closed.
\end{ex*}

\section{Integrable and half-integrable structures} 
\label{sec:integrable}
\noindent When the Nijenhuis tensor of an almost Hermitian manifold $N$ is
zero, the canonical torsion connection is the Bismut connection and the
$\G$T equations reduce to
\begin{equation*}
  \left\{
    \begin{aligned}
      d\psm&=\beta \psm,\\
      \omega d\omega&=\tfrac12\beta\oo.
    \end{aligned}
  \right.
\end{equation*}
The intrinsic torsion components (possibly) surviving are
\begin{equation*}
  W_3=\Omega,\qquad W_4=\tfrac{1}{4}\beta,\qquad W_5=-\tfrac12\beta,
\end{equation*}
the remaining data being encoded by $d\omega=\Omega+\tfrac12\beta\omega$.

The nasty expression \eqref{T} now becomes the much more tractable
\begin{equation*}
  T=-\hodge_6 d\omega-2\beta\hook\psm dt,
\end{equation*}
and a glance at \eqref{*T} shows

\begin{cor}\label{Tcoclosed}
  When $(N,J)$ is a complex manifold, the tensor $T$ is co-closed,\end{cor}

\noindent a fact implicitly noticed in~\cite{Friedrich-I:skew}.

\bigbreak
But our investigation focuses on six-dimensional nilmanifolds: they have
this striking property that makes all the above relations completely
trivial.  The holomorphic section $\Psi$ of the canonical bundle can be
expressed locally as product of (1,0)-forms $\alpha^i$ --- given an
$\SU$-structure such that $J$ is complex --- with $d\alpha^i\in{\mathcal
I}(\alpha^1,\dots ,\alpha^{i-1})$ \cite{Salamon:complex-nil}. Thus
$d\Psi=0$.  Conversely, if $(N,J,\psp)$ is a nilmanifold and $\psp,\psm$
are closed then $(d\alpha^i)^{0,2}=0$, hence
$d\Lambda^{1,0}\subseteq\Lambda^{2,0}\oplus\Lambda^{1,1}$.  Since
$\nabla\pspm=\mp\beta\otimes\pspm$ by hypothesis, the $1$-form $\beta$
vanishes identically.  What is more, the Lee form
$\vartheta=-Jd^*\omega=-\beta$ ($d^*=-\hodge_6d\hodge_6$ being the formal
adjoint of the exterior differential) is zero as well, so only the
$\W_3$-component of the intrinsic torsion remains:

\begin{prop}\label{balanced}
  An $\SU$-nilmanifold $(N,\psp,\omega)$ satisfying the $\G$T equations
  \eqref{g2t} and having an integrable almost complex structure is
  necessarily balanced.\qed
\end{prop}

A somehow simpler situation will crop up in \S\ref{SGT}, where
$(d\omega)^{2,1}_0$ will vanish too, and the manifold will turn out to
be K\"ahler.  So on products of the type $N\times S^1$ the integrability
hypothesis proves very restrictive and one problem is to understand the
case where the connection forms on $N$ and $S^1$ behave in a more entangled
way.  

\bigbreak
Having seen that $N_J=0$ forces the $\G$ Lee form to vanish and that the
$\G$T-structure is no longer of the general type, we consider another
distinguished setting which is not so restrictive, namely that of
half-integrability.

\begin{defn}
  {\rm\cite{Chiossi-S:SU3-G2}} An almost Hermitian $6$-manifold is
  \emph{half-integrable} if it possesses a reduction to $\SU$ for which
  both $\psp$ and $\oo$ are closed.
\end{defn}

Although this kind of structure has become popular in a certain part of the
physical literature with the name `half-flat', it might be preferable to
refer to such geometry with the term half-integrable, as
in~\cite{Apostolov-S:K-G2}.  When $\G$-manifolds are built from
$\SU$-structures, many features of the $\G$-structure are obtained by
properties in six dimensions, and the way an admissible $3$-form is built
motivated the above definition. This turns out to be a useful notion
especially in connection to Hitchin's evolution
equations~\cite{Hitchin:forms} which give metrics with holonomy $\G$ in
dimension~\( 7 \) starting from any such half-integrable structure in
dimension~\( 6 \)~\cite{Chiossi:PhD}.

In terms of $\SU$ classes half-integrability amounts to
\begin{equation*}
  W_1^+=0,\quad W_2^+=0,\quad W_4=0=W_5.
\end{equation*}
The hypothesis $\lambda\neq 0$ in equations \eqref{g2tclasses} is necessary
to avoid integrability issues, which we have already dealt with, so the
relevant non-vanishing derivatives are
\begin{equation*}
  d\psm=\tfrac12\lambda\oo,\qquad d\omega=\Omega-\tfrac{3}{4}\lambda\psp.
\end{equation*}
For the sake of completeness we also write down 
\begin{equation}
  \label{T*T}
  T=-\hodge_6 d\omega-\tfrac12\lambda\psm \qquad \hodge
  T=(d\omega+\tfrac12\lambda\psp)dt, 
\end{equation}
in accordance with the fact that $T=0$ implies the reduction of the
holonomy.

While leaving the discussion of the properties of $dT$ to Section~\ref{SGT}, we
are able now to weaken Corollary \ref{Tcoclosed} a little

\begin{cor}
  If $(N,\omega,\psp)$ is half-integrable then $T$ is co-closed.
\end{cor}

Notice that the Lee form $\vartheta$ of the $\SU$-structure vanishes and
the $\G$ Lee form~$\theta$ becomes locally exact, as required by
the Killing spinor equation mentioned in Remark~\ref{KS}.

\begin{rmk}
  A different way to simplify \eqref{g2t} is to assume the closure of
  $\psp$ only. But this annihilates $\beta$ (since
  $\beta\wedge\psp=0\then\beta=0$) and makes the
  $\W_1\oplus\W_2$-component vanish identically, so in the present
  context $d\oo=0$ follows from $d\psp=0$. In other words we land on
  half-integrability once again.
\end{rmk}

\section{Half-integrable nilpotent Lie algebras}
\label{HLNLAS}
\noindent The $\G$T equations combined with half-integrability form a powerful set of
constraints
\begin{equation}
  \label{startingeqns}
  \left\{
    \begin{aligned}
      d\psm&=\tfrac{\lambda}{2}\oo\\
      d\psp&=0
    \end{aligned}
  \right.
\end{equation} 
and the aim is to try to detect all half-integrable nilpotent Lie algebras
generating $\G$T-structures in seven dimensions in the described way. We
remind the reader that $\lambda\neq 0$ is the overall assumption from now
onwards.

Consider a $6$-dimensional nilpotent Lie algebra $\g$. Any such has a
nilpotent basis $(E^i)$, i.e.~a basis of $1$-forms such that
\begin{equation*}
  dE^i\in\Lambda^2V_{i-1},\quad \text{where}\quad V_j
  =\text{span}_\mathbb R\{E^1,\dots , E^{j-1}\},
\end{equation*}
and the spaces $V_j$ filter the Lie algebra: $0\subset V_1\subset\dots
\subset V_5\subset V_6=\g^*$. We shall indicate the basis $E^1,\dots ,E^6$
informally by $1,\dots ,6$ when no confusion arises.

\begin{rmk}
  \label{ker d}
  By adjusting the above spaces, it is always possible to take $ V_i$ to be
  the kernel of $d$, for some~$ i $, hence assume $V_2\subseteq\Ker d$.
  This small observation underpins the relation between the algebraic
  nilpotent filtration and the geometry under study.
\end{rmk}

In order to find all possible $\SU$-structures satisfying our equations
we need to determine a different orthonormal basis $e^1,\dots ,e^6$ of
$\g^*$ for which the K\"ahler form $\omega$ and the holomorphic volume
$\Psi$ are those of \eqref{stdsus}.

Let $U_j=V_j\cap JV_j$ be the subspaces corresponding to the $V_j$'s
invariant for the almost complex structure, so
\begin{equation*}
  \begin{array}{ccccccc}
    0&\neq & U_4 & \subseteq & U_5 & \subset & U_6\\
    & &\vrule height 1.3ex depth 0ex\mkern 2mu\cap & &\cap &
    & \vrule height 1.3ex depth 0ex \mkern 3mu \vrule height 1.3ex depth 0ex\\
    & & V_4 & \subset & V_5 & \subset & V_6. 
  \end{array}
\end{equation*}
Since $dU_j\subseteq\Lambda^2V_{j-1}$, a dimension count tells that $U_5$
has real dimension~$4$, and more importantly
\begin{equation*}
  2=\dim_\mathbb C U_4\iff JV_4=V_4,
\end{equation*}
that is $U_4$ is maximal if and only if $V_4$ is $J$-invariant. This will
be the first major watershed of the discussion. Since $\dim(U_5^\perp\cap
V_5)=1$ we also have $\dim_\mathbb C U_3,\dim_\mathbb C U_2\in\{0,1\}$ and
$\dim_\mathbb C U_1=0$, whence $U_2$ has (complex) dimension $1$ exactly if
$JV_2=V_2$.

\smallbreak
Let us fix the $\SU$-structure \eqref{stdsus} and concentrate only on the
hypothesis $V_4=U_4$, and prove

\begin{lemma}
  \label{lemma Jinv}
  For any nilpotent Lie algebra, under the assumption $V_4=JV_4$ there is
  no solution to $d\psm=\tfrac{\lambda}{2}\oo$, $d\psp=0$ (when
  $\lambda\neq 0$).
\end{lemma}

\begin{proof}
 One can assume
\begin{equation*}
  \1,\2,\3,\4\in U_4,\quad \5\in V_5\cap U_4^\perp,\quad \6\in JV_5\cap
  U_4^\perp.
\end{equation*}
Since $\g^*=V_5\oplus\lan\6\ran$, from $d\psm=\tfrac{\lambda}{2}\oo$ we induce
$d(e^{13}+e^{42})=-\lambda(e^{34}+e^{12})\5;$ but while the left-hand side of
the latter lives in $\Lambda^3V_4$, $\5\in V_5$ so the other side belongs to
$\Lambda^2V_4\wedge V_5$, and in particular $\lambda=0$. 
\end{proof}

\noindent {\it Examples. a)} The Iwasawa manifold once more provides a simple
instance. As $\dim_\mathbb R\Ker d\,|_{\Lambda^3}=15$, the space of complex
structures has dimension 12 (see~\cite{Ketsetzis-S:Iwasawa}), so the
closure of $\psp$ entails that $d\psm$ is proportional to $1234$, and the
$\G$T equations do not hold. It is no coincidence that it is precisely the
`$\lambda$ equation' that fails.

\smallbreak
{\it b)} More generally, for all Lie algebras with $dV_4=0$ something
similar happens. $V_4^\perp$ is $J$-invariant too, and there is a basis
such that $V_4=\text{span}\{e^1,\dots ,e^4\},
V_4^\perp=\text{span}\{e^5,e^6\}$. Let
\begin{equation*}
\alpha^1=\1+i\2,\quad \alpha^2=\3+i\4,\quad \alpha^3=\5+i\6
\end{equation*}
be a basis of complex $(1,0)$-forms determined by the almost complex
structure $J$, so we may write
\begin{equation*}
  \Psi=\alpha^1\wedge\alpha^2\wedge\alpha^3=\alpha^{123}.
\end{equation*}
Since $d\alpha^3 \in \Lambda^2 \lan \alpha^1, \alpha^2, \ol{\alpha^1},
\ol{\alpha^2} \ran$
\begin{equation*}
  d\Psi=\alpha^1\wedge\alpha^2\wedge
  d\alpha^3\in\Lambda^4V_4\otimes\mathbb C=\text{span}_\mathbb C\{e^{1234}\}.
\end{equation*}
But $\tfrac12\oo=e^{1234}+e^{1256}+e^{3456}$ does definitely not belong to
the latter, so equations \eqref{startingeqns} cannot be solved unless
$\lambda=0$.  A somehow more elegant proof of this fact will follow from
Lemma \ref{b_1}.  

\bigbreak
The possibility that $\dim_\mathbb C U_4=1<\dim_\mathbb C U_5$ involves
more thinking.  Let us change perspective for a moment and note that for an
invariant $\SU$ structure equations \eqref{startingeqns} necessarily have
$\lambda$ constant and so can be rewritten as the single equation
\begin{equation}\label{eq:d-Psi}
  d\Psi=\tilde{\lambda}\oo,
\end{equation}
with $\tilde{\lambda}=i\lambda/2$. One may now regard this as an equation for a
$U(3)$-structure and allow $\tilde\lambda$ to be a complex number.

The symbol $\alpha^{\bar{\imath}}$ will denote $\overline{\alpha^\imath}$ for
$i=1,2,3$.  Since $\Psi=\alpha^{123}$ and
\begin{equation*}
  -2\omega^2 = \alpha^{1\bar 12\bar 2} + \alpha^{1\bar 13\bar 3}
  + \alpha^{2\bar 23\bar 3}, 
\end{equation*}
equation~\eqref{eq:d-Psi} becomes
\begin{equation*} 
(d\alpha^1)\alpha^{23} + (d\alpha^2)\alpha^{31} + (d\alpha^3)\alpha^{12}
  = -\tfrac{1}{2}\tilde\lambda
  (\alpha^{1\bar 12\bar 2} + \alpha^{1\bar 13\bar 3} + \alpha^{2\bar 23\bar 3}),
\end{equation*}
which implies immediately
\begin{equation}
  \label{eq:daj-02}
  \begin{aligned}
    (d\alpha^1)^{0,2} &= \tfrac12 \tilde\lambda \alpha^{\bar 2\bar 3},\\
    (d\alpha^2)^{0,2} &= \tfrac12 \tilde\lambda \alpha^{\bar 3\bar 1},\\
    (d\alpha^3)^{0,2} &= \tfrac12 \tilde\lambda \alpha^{\bar 1\bar 2}.
  \end{aligned}
\end{equation}
Thus a choice of e.g.~$\alpha^1$ determines the span of $ \alpha^{\bar 2} $
and $ \alpha^{\bar 3} $ in $ \Lambda^{0,1}$.  Conjugating, one gets the
span of $ \alpha^2 $ and $ \alpha^3 $ in~$ \Lambda^{1,0} $, and hence the
full orthogonal complement of $
\text{span}\{\alpha^1,\alpha^{\bar 1}\}$.

As $ U_4, V_3\subset V_4 $ and these have real dimensions $ 2 $, $3$ and~$
4 $, respectively, we have
\begin{equation*}
  \dim U_4 \cap V_3 \geqslant 1.
\end{equation*}
Similarly,
\begin{gather*}
  \dim U_5 \cap U_4^\perp \cap V_4 \geqslant 1,\\
  \dim U_5^\perp \cap V_5 \geqslant 1.
\end{gather*}
We may now take 
\begin{alignat*}{3}
  \alpha^1 &= e^1+ie^2 \in U_4 &\quad&\text{with}\quad& e^1 &\in V_3,\\
  \alpha^2 &= e^3+ie^4 \in U_5\cap U_4^\perp &&\text{with}& e^3 &\in V_4,\\
  \alpha^3 &= e^5+ie^6 \in U_5^\perp &&\text{with}& e^5 &\in V_5.
\end{alignat*}
The first of equations~\eqref{eq:daj-02} implies that $ d\alpha^1\ne0 $, so
one of $de^1$, $de^2$ is non-zero and therefore

\begin{lemma}
  \label{b_1}
  The space of closed $1$-forms has real dimension at most three, i.e.~$
  b_1\leqslant 3 $.
  \qed
\end{lemma}

Notice by the way that the first Betti number of any nilmanifold is always
strictly greater than one~\cite{Dixmier:nil-cohomology}, which reduces our
investigations to the nilpotent Lie algebras with $b_1$ equal to either $2$
or $3$. It becomes clear how crucial Remark~\ref{ker d} is in the
description, for it now states that the kernel of the exterior derivative
is $V_2$ or possibly $V_3$.

\medbreak Although the argument of Lemma \ref{lemma Jinv} can be adapted to
fit Lie algebras with many closed forms, a straightforward consequence of
the previous result is that

\begin{cor}
  For any six-dimensional nilpotent Lie algebra with $b_1\geq 4$ the
  equations $d\psm=\tfrac{\lambda}{2}\oo,d\psp=0$ have no solution, for
  $\lambda$ non-zero.  \qed
\end{cor}

\noindent The complex filtration $\{U_\alpha\}$ is rigid in the sense the
size of the spaces $U_\alpha$ with odd index is fixed, since $\dim_\mathbb
C U_1=0$, $\dim_\mathbb C U_5=2$, and we claim that
\begin{lemma}
  $ U_3=\{0\} $ when $U_4$ is not maximal.
\end{lemma}
\begin{proof}
 If $ U_3$ were non-trivial then $ e^1,e^2 \in V_3 $ and we might
take $ e^1 \in V_2 $. Then $ d\alpha^1 \in\Lambda^2V_2 $ and $ e^1d\alpha^1
= 0 $.  But
\begin{equation*}
(e^1d\alpha^1)^{0,3} = (e^1)^{0,1}(d\alpha^1)^{0,2} = \tfrac 14\tilde \lambda
\alpha^{\bar1\bar2\bar3} \ne 0, \end{equation*} yielding a contradiction.
\end{proof} With $
\dim_\mathbb C U_3 =0$, we now have
\begin{gather*}
  e^2 \in V_4\setminus V_3,\\
  e^4 \in V_5\setminus V_4,\\
  e^6 \in V_6\setminus V_5.
\end{gather*}
As $ \{e^1,e^2,e^3,e^4,e^5 \} $ is a real orthonormal basis for $ V_5 $, we
conclude
\begin{equation*}   e^6 \in V_5^\perp. \end{equation*}
In real components equation~\eqref{eq:d-Psi} is 
\begin{equation*}
  d\bigl((e^1+ie^2)(e^3+ie^4)(e^5+ie^6)\bigr) =
  2\tilde\lambda(e^{1234}+e^{1256}+e^{3456}),
\end{equation*}
whose $ e^6 $-component gives
\begin{equation}\label{eq:6}
  d\bigl( (e^1+ie^2)(e^3+ie^4)\bigr) = 2\tilde\lambda (e^{12}+e^{34})e^5.
\end{equation}
As $ e^4 \in V_5\setminus V_4 $ and $ e^5 \in V_5 $ we may write
\begin{equation*}
  e^5 = qe^4 + \xi,
\end{equation*}
for some $ q\in\mathbb R $ and $ \xi\in V_4 $. Note that $\xi$ completes
$\{e^1,e^2,e^3\} $ to an orthogonal basis for~$ V_4 $.
Equation~\eqref{eq:6} is now
\begin{equation*}
  d\bigl((e^1+ie^2)(e^3+ie^4)\bigr) = 2\tilde\lambda
  \bigl((qe^{12}-e^3\xi)e^4+e^{12}\xi\bigr)
\end{equation*}
and extracting the $ e^4 $-component we have
\begin{equation*}
  - d(e^1+ie^2) = 2\tilde\lambda(qe^{12}-e^3\xi).
\end{equation*}
However $ d(e^1+ie^2)\in \Lambda^2V_3\otimes\mathbb C $ and as $ \dim V_3=3
$ this two-form is necessarily decomposable.  This is true of the
right-hand side only if $ q=0 $.  Therefore we conclude firstly that
\begin{gather*}
  e^5 \in V_4,\\
  e^4 \in V_5\cap V_4^\perp.
\end{gather*}
Our equation now reads $ d(e^1+ie^2)=2\tilde\lambda e^{35} $.  Thus some
real linear combination of $ e^1 $ and $ e^2 $ is closed.

\begin{lemma}\label{d2}
  \begin{equation*}
    de^1=0,\qquad de^2 = \lambda e^{35}.
  \end{equation*}
\end{lemma}

\begin{proof}
  This follows directly from the fact that one of the $V_i$'s is the kernel
  of $d$ and $ \lambda = -i2\tilde\lambda$ is real.
\end{proof}

Because $ \{e^1,e^2,e^3,e^5\} $ is a basis for $ V_4 $, and $ e^2\in
V_4\setminus V_3 $, we may write
\begin{gather*}
  e^3 = c_1e^2 + \eta_1,\\
  e^5 = c_2e^2 + \eta_2,
\end{gather*}
with $ c_i\in\mathbb R $ and $ \eta_j $ linearly independent vectors
of $ V_3 $.  By Lemma~\eqref{d2}, 
\begin{equation*}
  \Lambda^2V_3 \ni e^{35} = (c_2\eta_1-c_1\eta_2)e^2 + \eta_1\eta_2.
\end{equation*}
Thus $ c_2\eta_1-c_1\eta_2=0 $, which gives $ c_1=0=c_2 $, and hence
\begin{gather*}
  e^3, e^5 \in V_3,\\
  e^2 \in V_4\cap V_3^\perp. 
\end{gather*}
In particular, $ V_3 $ has orthonormal basis $ \{e^1,e^3,e^5\} $ and $
V_3^\perp $ has orthonormal basis $ \{e^2,e^4,e^6\} $, so $
V_3^\perp=JV_3$:

\begin{cor}
  \begin{equation*}
    \g^*=V_3\oplus JV_3=\lan\1,\3,\5\ran\oplus\lan\2,\4,\6\ran.
  \end{equation*}\qed
\end{cor}
 \noindent With this, we are able to pin down the bases for almost every
space in the filtration, with the exception of the first two. Indeed
\begin{align*}
  V_3&=\lan\1,\3,\5\ran,\\
  V_4&=V_3\oplus\lan\2\ran,\\
  V_5&=V_4\oplus\lan\4\ran. 
\end{align*}
In the light of Remark~\ref{ker d}, we could have chosen the filtration so that
to fix $V_1, V_2$ as well. As we shall see at the end of the section, it
would have been possible, in theory, to decree $V_1$ to be the span of $\1$, and
$V_2$ to be either $\lan\1,\3\ran$ or $\lan\1,\5\ran$. The computation will
indeed show that $\3$ and $\5$ play interchangeable r\^oles. Apparently
though not much is gained from this fact, a reason why we prefer the more
general point of view.

We now inspect the space $V_3$, with regard to the fact that
$dV_3\subseteq\Lambda^2V_2$. Since $d\3\in\Lambda^2V_2$ either $\3$ is
closed or $b_1=2$ (recall that $\1\in\Ker d$ anyway), so $\1\in V_2$.
Whichever the case we may write
\begin{equation}\label{1d3}
  \1 d\3=0.
\end{equation}
A similar argument holds for $d\5\in\Lambda^2V_2$, whence 
\begin{equation}\label{1d5}
  \1 d\5=0.
\end{equation}
These two relations will turn out useful later on.

Let us look at the real and imaginary parts of \eqref{eq:6}, now giving the
reassuring
\begin{gather*}
  d(e^{14}+e^{23})=0,\\
  d(e^{13}+e^{42})=\lambda(e^{12}+e^{34})\5,
\end{gather*}
or, by means of Lemma \ref{d2}, simply
\begin{gather}
  \1 d\4+\2 d\3=0\\ 
  \2 d\4=\lambda e^{125}.\label{x}
\end{gather}

Besides, from \eqref{x} we may write $d\4=-\lambda e^{15}+\2\wedge V_3$.
Explicitly setting $d\4=-\lambda e^{15}+\2\wedge(v\1+u\3+t\5)$ gives
\begin{equation}\label{pre nu1}
\1(u\3+t\5)=d\3. \end{equation}

Since the codimension of $V_2$ in $V_3$ is one, some linear combination of
$\3$ and $\5$ is in $V_2$, hence closed. In fact, fix a non-zero element
$\nu$ in $\Lambda^2V_2=\mathbb C$, so that
\begin{equation}
  \label{nu}
  d\3=x\nu,\quad d\5=y\nu
\end{equation}
for some $y,x\in\mathbb R$. Then the closure of $e^{35}$ makes
$(x\5-y\3)\nu$ vanish, whence
\begin{equation*}
  x\5-y\3\in V_2.
\end{equation*}

The equation $\omega d\omega=0$ is a direct consequence of the
first in \eqref{startingeqns}, so we should not expect to gain new 
information by manipulating it. Nonetheless, it quickly furnishes 
substantial results
\begin{equation*}
  d(e^{12}+e^{34})\5+(e^{12}+e^{34})d\5=0,
\end{equation*}
in which we separate $\4$-terms from the remaining ones above and obtain
\begin{equation}
  e^{35}d\4=0.\label{3'perp}
\end{equation}
Then \eqref{pre nu1} rephrases as 
\begin{equation}
  \1(u\3+t\5)=x\nu\label{nu1}
\end{equation}
whilst \eqref{3'perp} implies $v=0$, so
\begin{equation*}
  d\4=-\lambda e^{15}+ue^{23}+te^{25}.
\end{equation*}

More is recovered by considering the $V_5$-component of \eqref{eq:d-Psi},
whose real and imaginary parts now become
\begin{equation}\label{V5}
  \begin{gathered}
    -e^{24}d\5 = (e^{14}+e^{23})d\6, \\
    -e^{15}d\4 +(e^{13}+e^{42})d\6 = \lambda e^{1234}.
  \end{gathered}
\end{equation}
Defining
\begin{equation*}
  d\6=\alpha\wedge \4+\zeta, \qquad\text{with}\quad \zeta\in\Lambda^2V_4,
  \quad \alpha\in V_4,
\end{equation*}
the $(\4)^\perp$- and $\4$-terms in \eqref{V5} yield
\begin{equation} \label{newtie}
  \left\{
    \begin{aligned}
      e^{13}\zeta & =ue^{1235},\quad & e^{13}\alpha-\2\zeta & =\lambda e^{123},\\
      e^{23}\zeta & =0,              & \1\zeta+e^{23}\alpha & =-y\2\nu.
   \end{aligned}
  \right.
\end{equation}
At this very moment it seems difficult to get much out of these relations.
Still we are able to say that the forms $\zeta$ and $\alpha$ are of the kind
\begin{gather*}
  \zeta=z_2e^{12}+z_3e^{13}+z_1e^{23}-ue^{25}\\
  \alpha=(z_3-\lambda)\2+a_1\1+a_3\3
\end{gather*}
with $z_i,a_j$ real constants. Moreover, the last equation of
\eqref{newtie} tells that
\begin{equation}\label{y}
  y\nu=-(z_1+a_1)e^{13}+ue^{15}.
\end{equation}
The constraint $d(d\6)=0$ forces $z_2=a_3$, whence
\begin{equation}
  \label{dd6}
  d\6=a_1e^{14}+z_1e^{23}-a_3(e^{12}-e^{34})+z_3(e^{13}-e^{42})+\lambda
  e^{42}-ue^{25}
\end{equation}
plus other complicated relations whose full exploitation we will delay. 

Given equations \eqref{1d3} and \eqref{1d5}, we are obliged to separate the
discussion. It is preferable to handle mutually exclusive cases, so we
shall assume that $b_1=2$ in the first two cases, allowing three
independent closed $1$-forms in the third only.
\begin{description}
\item[Case 1] is that in which we suppose that $\3$ is not closed. 
\item[Case 2] instead has the derivative of $\3$ zero but $\5$ non-closed.
\item[Case 3] deals with both $\3,\5$ being closed (so $b_1=3$).
\end{description}
Note that equations \eqref{nu1} and \eqref{y} do not permit $\3,\5$ to be
both non-closed. Since $\nu=\1(x\5-y\3)$ in fact, if $x,y$ were non-zero we
would have $u=-xy=xy$, i.e.~$\5\in V_2\subseteq\Ker d$.

\subsection*{Case 1}
\noindent Let $k$ be the non-zero real number such that $\Lambda^2V_2\ni
d\3=ke^{15}$.  The $2$-form $\nu$ of \eqref{nu} is merely proportional to
$e^{15}$, so from \eqref{nu1} we infer that $u=0,t=k$. In addition, by
\eqref{y} we also have $a_1=-z_1$.

The new structure relations are
\begin{equation*}
  d\3=ke^{15},\quad d\4=-\lambda e^{15}+ke^{25},\quad d\5=0,
\end{equation*}
so the string $(de^i)$ only lacks the explicit determination of the last
entry. Instead of chasing every little piece of data around, we use
\eqref{newtie}.  The demand $d(d\6)=0$ reduces the number of coefficients
in \eqref{dd6}, in fact $z_3=\lambda,\ a_3=0$ and $a_1=z_1 (=-a_1$), thus
$d\6=\lambda e^{13}$ and finally the Lie algebra structure is:
\begin{equation}\label{d3neq0}
  (0,\lambda e^{35},ke^{15},-\lambda e^{15}+ke^{25},0,\lambda e^{13}),
\end{equation} 
with $\lambda\neq 0\neq k$.

\smallbreak
The next task is to spot this Lie algebra within the list
of~\cite{Salamon:complex-nil}, entries on which will be numbered in bold
starting from the top, so for example \textbf{28} indicates the Lie algebra
of the Iwasawa manifold, \textbf{34} that of a torus.

Let us compute the Betti numbers of~\eqref{d3neq0}: while $b_1$ is clearly
$2$, we have that $\dim\Ker d\,|_{\Lambda^2}=8$, so $b_2=4$. The only Lie
algebra with such cohomological data and decomposable exact $2$-forms is
number \textbf{6}, in other words
\begin{equation*}
  (0,0,12,13,23,14).
\end{equation*}
This could have been checked by modifying the basis in \eqref{d3neq0}. In
fact, redefining $\4$ as $\lambda\3+k\4$ gives $(0,\lambda
e^{35},ke^{15},k^2e^{25},0,\lambda e^{13})$. The successive swaps
$\2\leftrightarrow\5,\6\leftrightarrow\4,\1\leftrightarrow\2,
\4\leftrightarrow\5$ followed by a flip in the signs of $\1,\6$ lead to
\textbf{6}, with suitable rescaling.

\subsection*{Case 2}
\noindent This case will take slightly longer, due
to minor complications. Given that $V_2=\text{span}\{\1,\3\}$, equation
\eqref{nu1} provides a neat expression for the derivative of $\4$, namely
\begin{equation*}
  d\4=-\lambda e^{15}.
\end{equation*}
In addition \eqref{y} determines 
\begin{equation*}
  d\5=-(z_1+a_1)e^{13},
\end{equation*}
with the assumption that $z_1+a_1\neq 0$.  When asking $d\6$ to be closed,
we find that $z_3=\lambda$, so we may write the Lie algebra as
\begin{equation*}
  \bigl(0,\lambda e^{35},0,-\lambda
  e^{15},ze^{13},\tfrac{a_1+z}{2}(e^{14}+e^{23})
  +\tfrac{a_1-z}{2}(e^{14}-e^{23})-a_3(e^{12}-e^{34})+\lambda
  e^{13}\bigr), 
\end{equation*}
where $z$ stands for $z_1$. 

The fact that now $d\4$ lives in $\Lambda^2V_3$ renders the filtration
$\{V_i\}$ more flexible to further adjustments. The $\SU$-structure fixes
the span of $\alpha^3$, hence $\5,\6$ are confined by this diagram
\begin{equation*}
  \begin{array}{ccccccc}
    V_2 & \subset & V_3 & \subset & V_5 & \subset & V_6 \\
    \1,\3 & & \5 & & \2,\4 & & \6.
  \end{array}
\end{equation*}
Furthermore, if $\mathbb D$ denotes the real space
$\lan\1,\ldots,\4\ran$, then
\begin{equation*}
  \Lambda^2{\mathbb D}=\Lambda^2_+\oplus\Lambda^2_-,
\end{equation*}
where $\Lambda^2_-$ is generated by the anti-self-dual forms
$e^{14}-e^{23},e^{12}-e^{34},e^{13}-e^{42}$, and similarly for the self-dual
part.  Since the decomposition
\begin{equation*}
  \mathbb D=\lan\1,\3\ran \oplus \lan\2,\4\ran
\end{equation*}
has to be preserved, it is possible to act by $\Lie{SO}(2)\subset
\Lie{SU}(2)$ on the $2$-plane $\lan e^{14}-e^{23},e^{12}-e^{34}\ran\subset
\Lambda^2_-$ to eliminate the $a_3$ term above.  More explicitly, the
required transformation on $\mathbb D=\lan\alpha^1,\alpha^2\ran_{\mathbb R}$ is given by
\begin{equation*}
  \alpha^1\mapsto\alpha^{\tilde 1} = c\alpha^1+s\alpha^2 \qquad 
  \alpha^2\mapsto\alpha^{\tilde 2} = -s\alpha^1+c\alpha^2,
\end{equation*}
for some rotating constants $c=\cos \sigma$, $s=\sin \sigma$.  The common
coefficient of $e^{\tilde 1\tilde 2},e^{\tilde 3\tilde 4}$
\begin{equation*}
  sc(z-a_1)-a_3(s^2-c^2)
\end{equation*}
is killed off by a suitable choice of angle, given by $\tan
2\sigma=(a_1-z)/2a_3$. 
That allows to simplify the general expression of $d\6$, and finally the
Lie algebra looks like
\begin{equation}\label{S6-S7-S8}
  (0,\lambda e^{35},0,-\lambda e^{15},ze^{13},a_1e^{14}-ze^{23}+\lambda e^{13}).
\end{equation}

Now, we need only determine the second Betti number: whilst there are six
obvious closed $2$-forms, $de^{16}=-ze^{123}$ and $de^{36}=a_1e^{134}$ are
proportional to $de^{25},de^{45}\in \im d\,|_{\Lambda^2}$, irrespective of
the coefficients. Therefore $b_2=4$ and \eqref{S6-S7-S8} possibly
identifies with \textbf 6, \textbf 7, \textbf 8 of Salamon's list.

For instance, $a_1=0$ gives \textbf 6, in which every exact $2$-form is
decomposable. It is easy to find an explicit isomorphism in this case.
First, change $\6$ to $z\6+\lambda\5$ to get rid of the redundant last term
in $d\6$, then swap $\2,\3$ in order to have $(0,0,\lambda e^{25},-\lambda
e^{15},-ze^{12},-ze^{23})$. Now, $d\4\wedge d\6=z\lambda e^{1235}$ is the
only non-zero wedge product of exact $2$-forms (while in \textbf 6
$d5\wedge d6=1234$), suggesting to swap $\4$ and $\5$.  Then it is a
question of exchanging $\3,\4$ and $\1,\2$, reversing the orientation of
$\5,\6$ and setting $z=1$.

\medbreak
As far as $(0,0,12,13,23,14\pm 25)$ are concerned, their characteristic
series agree with those of \eqref{S6-S7-S8} regardless of coefficients: the
central descending series' dimensions being $6,4,3,1,0$, the derived
series' $6,4,0$ and the upper central series' ones $1,3,4,6$. But this is
not enough to distinguish the sign, and a more subtle argument has to be
given. By analysing \textbf 7, \textbf 8 we can say a few things about the
generic $\SU$ basis, namely that $\Ker d=\lan 1, 2\ran$, and since
$d3\wedge\g^*=0$, the form $3$ is specified up to a choice of $1, 2$. In
addition, $d4, d5$ determine $4, 5$ (up to $1, 2$), and have a common
factor $3$, hence $3$ is now completely determined. As for $d6$, it just belongs
to $\lan 1, 2\ran\otimes\lan 4, 5\ran$.

Furthermore, the only non-vanishing relations among exact $2$-forms are
$d4\wedge d6=\mp 1235,d5\wedge d6=1234$.  Let $x4+y5$ be the candidate
for fifth basis element. Then
\begin{equation*}
  d(x4+y5)\wedge d6=\mp x1235+y1234=123(\mp x5+y4)
\end{equation*}
from which $\mp x5+y4$ is the new $4$, so
\begin{equation*}
  d(\mp x5+y4)\wedge d6=\mp 123(x4+y5).
\end{equation*}
In \eqref{S6-S7-S8} instead, after renaming $\6=(z+a_1)\6+\lambda\5$, we
have
\begin{align*}
  d\2\wedge d\6 &= \lambda e^{13}(-a_1\4)\5,\\
  d\4\wedge d\6 &= \lambda\1(-z\2)e^{35};
\end{align*}
so apparently $\lan 1,2\ran$ corresponds to $\lan \1,\3\ran$, $3$ should be
$\5$ and $\lan 4,5\ran$ is $\lan \2,\4\ran$.  Comparing
\begin{gather*}
  d(-a_1\4)\wedge d\6 = \lambda a_1z e^{1235}\qquad\text{and}\\
  d(-z\2)\wedge d\6 = \lambda a_1z e^{1345}
\end{gather*}
with the corresponding $d4\wedge d6=\mp 1235, d5\wedge d6=1234$, relative
to \textbf 8, \textbf 7, shows that $\lambda a_1z$ can be both positive or
negative, so the algebra type is detected by the sign of this particular
coefficient.

\begin{cor}
  All nilpotent Lie algebras with Betti numbers $b_1=2,b_2=4$
  \begin{gather*}
    (0,0,12,13,23,14)\\
    (0,0,12,13,23,14+25)\\
    (0,0,12,13,23,14-25)
  \end{gather*}
  bear an $\SU$-structure fulfilling \eqref{startingeqns}.
\end{cor}

It turns out quite instructive to provide a different proof of the result
by exhibiting a change of basis. The detailed description given for \textbf
7, \textbf 8 suggests that one should first exchange $\2$ with $\3$, then
$\5$ with $\3$. Then redefine $\6$ as the linear combination
$(z+a_1)\6+\lambda\3$ in order to have $(z+a_1)d\6=a_1e^{14}-ze^{25}$. The
algebra $\bigl(0,0,-(z+a_1)e^{12},-\lambda e^{13}, \lambda e^{23},
(z+a_1)a_1e^{14}-z(a_1+z)e^{25}\bigr)$ is almost of the desired form.
Forgiving the abuse of notation, a diagonal transformation
\begin{equation}\label{diagonal}
  \1 \mapsto a\1,\ \2 \mapsto b\2,\ \3 \mapsto
  c\3,\ \4 \mapsto f\4,\ \5 \mapsto g\5,\ \6 \mapsto h\6 
\end{equation}
allows enough freedom to detect the correct coefficients: for instance
$d\3$ becomes $-(z+a_1)(c/ab)e^{12}$ etc., and normalisation entails
\begin{gather*}
  -(z+a_1)c=ab,\quad -\lambda f=ac,\quad \lambda g=bc\\
  (z+a_1)a_1h=af,\qquad -z(a_1+z)h=\pm bg.
\end{gather*}
Thus $c,f,g,h$ are determined by $a,b$ and $z=\pm a_1b^2/a^2$.  Taking
$a=1$ it is possible to assign $a_1$, $z$, find $b$ accordingly, and so on
with the remaining numbers.

\subsection*{Case 3}
\noindent Since now each form in $V_3$ is closed, \eqref{pre nu1} 
determines $d\4=-\lambda e^{15}$, and similarly
\eqref{y} yields $z_1=-a_1$ in \eqref{dd6}. Recall that $a_3$ can be taken
to be zero. Imposing $d\6$ to be closed further gives $z_3=\lambda$, so one
eventually arrives at
\begin{equation}\label{d3=0=d5}
  \bigl(0,\lambda e^{35},0,-\lambda e^{15},0,a_1(e^{14}-e^{23}) 
  +\lambda e^{13}\bigr).
\end{equation}
The cohomology depends on the value of $a_1$ as suspected. Indeed, the
inspection of the differentials $de^{ij}$ says that $\dim\Ker
d\,|_{\Lambda^2}=7$.  Thus generically $b_2=5$ whilst $b_1=3$, giving a
large range to choose from, in principle. But when $a_1$ vanishes, $b_2$
raises to 8 and the Lie algebra is
\begin{equation*}
  (0,0,0,12,13,23)
\end{equation*}
in disguise. This corresponds to the example of a non-locally conformally
flat $\SU$-structure carrying parallel torsion given
in~\cite{Ivanov-I:instantons}.

\medbreak
If $a_1$ is non-zero instead, two non-isomorphic real Lie algebras crop up
\begin{cor}
  Of all nilpotent Lie algebras with $b_1=3,b_2=5$ only
\begin{equation*}(0,0,0,12,23,14\pm 35)\end{equation*}
possess half-integrable structures satisfying \eqref{startingeqns}.
\end{cor}

\begin{proof}
Assume $b_2=5$. First of all~\eqref{d3=0=d5} is $3$-step for all choices
of $a_1$, which rules out some possibilities. Secondly the dimensions of
the central descending series are $6,3,1,0$, so we are only left with
\textbf{14}, \textbf{15}, \textbf{16} (the other series not distinguishing
further).

Now we start playing around with bases: swapping $\2,\5$ allows to write
\begin{equation}
  \label{15-16}
  \bigl(0,0,0,-\lambda e^{12},-\lambda e^{23},\lambda e^{13} 
  +a_1(e^{14}+e^{25})\bigr).
\end{equation}
Consider the following relations relative to the Lie algebras on the
right
\begin{equation*}
  \textbf{14}
  \left\{
    \begin{aligned}
      &d4\wedge d6=1235,\\
      &d5\wedge d6=0,
    \end{aligned}
  \right. \qquad
  \textbf{15}, \textbf{16}
  \left\{
    \begin{aligned}
      &d4\wedge d6=\pm 1235,\\
      &d5\wedge d6=1234.
    \end{aligned}
  \right.
\end{equation*} 
These suggest to redefine $\4$ by a linear combination $x\4+y\5$. Then
$d(x\4+y\5)\wedge d\6=\lambda( -ya_1)e^{1234}+\lambda(-xa_1)e^{1235}$
imply $x=1,y=0$ hence that $\4$ should remain unchanged, but $\5$ itself
should be rescaled by $-a_1\5$. The product $d(-a_1\5 )\wedge d\6=\lambda
a_1^2 e^{1234} $ tells that \eqref{d3=0=d5} corresponds to
\begin{alignat*}{3}
  \textbf{14}&&\quad&\text{if and only if}&\quad&a_1 = 0,\\
  \textbf{15}, \textbf{16}&&&\qquad\text{if}&&a_1^2=1.
\end{alignat*}
While the vanishing of $a_1$ is excluded by assumption, $a_1=\pm 1$ will
transform \eqref{d3=0=d5} into either \textbf{16} or \textbf{15}.
\end{proof}

\noindent Explicitly, replacing $\4$ in $d\6$ with $\4+a\1+b\3$ gives 
$d\6=(\lambda+a_1b)e^{13}+a_1(e^{14}+e^{35})$, 
simplified to $d\6=a_1(e^{14}+e^{35})$ when $b=-\lambda/a_1$. Whilst
\textbf{15} crops up once we set $1=a_1$, we apply the transformation
\eqref{diagonal} to \eqref{15-16}, and then normalise coefficients
\begin{equation*}
  -\lambda f=ab,\quad -\lambda g=bc,\quad
  a_1h=af,\quad -a_1h=cg.
\end{equation*}
Choosing $a=1,b=\lambda,c=1$, tantamount as reversing the orientations of
$\4,\5$ and rescaling $\2,\6$, eventually produces \textbf{16}.

\bigbreak
We are now able to collect the outcome of all previous Corollaries together 

\begin{thm}\label{output}
  Up to isomorphism, there are precisely 6 real six-dimensional
  half-integrable nilpotent Lie algebras satisfying the $\G$T equations
  \begin{gather*}
    (0,0,12,13,23,14)\\
    (0,0,12,13,23,14-25)\\
    (0,0,12,13,23,14+25)\\
    (0,0,0,12,23,14+35)\\
    (0,0,0,12,23,14-35)\\
    (0,0,0,12,13,23)
  \end{gather*}
  with non-zero Lee form $\theta$.\qed
\end{thm}

\noindent Observe that these involve only 4 different complex types. As the real
parameter $\lambda$ is fixed, each of \eqref{d3neq0}, \eqref{S6-S7-S8},
\eqref{d3=0=d5} is defined in terms of effective coefficients, varying
which changes the $\SU$ reduction and possibly the Lie algebra. Therefore
the complex types
\begin{equation*}
  (0,0,12,13,23,14)\qquad (0,0,0,12,23,14\pm 35)\qquad (0,0,0,12,13,23)
\end{equation*}
admit a one-parameter family of half-integrable $\SU$-structures, whilst 
\begin{equation}\label{proto}
  (0,0,12,13,23,14\pm 25)
\end{equation}
possesses a real $2$-dimensional family of such.

\smallbreak
Finally, note that the latter twin algebras, which have highest step-length
and smallest Betti numbers, can degenerate to all others by an appropriate
contraction limit~\cite{Gibbons-LPS:domain-walls}. As a typical example, we
show the contracting procedure by means of which $(0,0,12,13,23,14)$ is
attained as limit of~\eqref{proto}. Introduce the real parameter $t$ and
choose the following basis for $\Lambda^1$ depending on $t$:
\begin{equation*}
  t^{-1}\1,\quad t\2,\quad \3,\quad t^{-1}\4,\quad t\5,\quad t^{-2}\6.
\end{equation*}
This leaves all differentials $de^i$ unchanged except for $d\6=e^{14}\pm
t^{-4}e^{25}$, which becomes $e^{14}$ when $t\lto 0$. This is same as setting
$a_1$ equal to zero in~\eqref{S6-S7-S8}. The technique is
also designed to increase Betti numbers, and produces the other Lie
algebras of Theorem~\ref{output} in a natural way.

\section{The derivative of the torsion}
\label{SGT}
\noindent The behaviour of $dT$ is examined here on a case-by-case
basis. The Lie algebras \eqref{d3neq0} and~\eqref{S6-S7-S8} produce similar
expressions
\begin{gather*}
  dT=\tfrac 32\lambda^2\oo-2k^2e^{1256},\\
  dT=\tfrac 32\lambda^2\oo-(a_1^2+2z^2)e^{1234},
\end{gather*}
both of which have a $\oo$-component plus a non-vanishing extra term.  
For the algebras
encountered in Case~3 no additional coefficient appears, hence the
derivative of $T$ has the simplest possible form
\begin{equation*}
  dT=\tfrac 32\lambda^2\oo.
\end{equation*}

\begin{lemma}\label{eform}
  The K\"ahler form of the Lie algebras
  \begin{equation*}
    \bigl(0,\lambda e^{35},0,-\lambda e^{15},0,a_1(e^{14}-e^{23})
    +\lambda e^{13}\bigr)
  \end{equation*}
  is an eigenform of the Laplace-Beltrami operator $\Delta$.
\end{lemma}

\begin{proof}
  Since $d\hodge_6 d\omega=-\tfrac 32\lambda^2\oo$ we have
  $d^*d\omega=-\hodge_6 d\hodge_6 d\omega=3\lambda^2\omega$. But the
  squared K\"ahler form $\oo$ is closed and thus $d^*\omega=-\hodge_6
  d(\tfrac12\oo)=0$, so the Laplacian
  \begin{equation*}
    \Delta=d^*d+dd^*
  \end{equation*}
  acts on $\omega$ just as the operator $d^*d$.
\end{proof}

\noindent Moreover, it is always true that

\begin{cor}
  The exterior derivative of the torsion $T$ is a form of type $2,2$ with
  respect to $J$.
\end{cor}

\noindent The dependence on the choice of one particular $J$
vaguely reminds of the analogous situation in the quaternionic K\"ahler
case, where a strong QKT structure is purposely defined by
$dT\in\Lambda^{2,2}$, this time though with respect to the whole sphere of
almost complex structures~\cite{Howe-OP:QKT}.

Recall that in terms of the fundamental $\G$-representation $V_7$
\begin{align*}
  T\in\Lambda^3 & \iso  \mathbb R\oplus V_7\oplus\sym^2_0V_7\\
  dT\in\Lambda^4 & \iso  \mathbb R\oplus V_7\oplus\sym^2_0V_7.
\end{align*}
Since $T$ is a tensor attached to seven-dimensional geometry, a better
counterpart is given by
\begin{prop}
  \begin{equation*}
    dT\in\mathbb R\oplus\sym^2_0V_7\subset\Lambda^4,
  \end{equation*}
\end{prop}

\noindent improving the previous Corollary.

\begin{proof}
 It is well-known that
\begin{align*}
  V_7         & =\rcomp{\Lambda^{1,0}} \oplus \mathbb R,\\
  \sym^2_0V_7 & =\rcomp{\sym^{2,0}} \oplus \rreal{\Lambda^{1,1}_0} \oplus
  \rcomp{\Lambda^{1,0}} \oplus \mathbb R 
\end{align*}
as irreducible $\SU$ modules. Take a tangent vector $X$ to $M$, so
$X=\tfrac{\partial}{\partial t}+\hat X$ with $\hat X\in TN$.  Then
\begin{equation*}
  \begin{split}
    (X\hook\hodge\varphi)\wedge\oo
    & =(\tfrac 12\hat X\hook\oo-\psm+\hat X\hook\psm dt)\oo\\
    & =(\tfrac 12\hat X\hook\oo)\oo-\psm\oo+(\hat X\hook\psm dt)\oo=0.
  \end{split}
\end{equation*}
The first two terms are $7$-forms on a six-dimensional manifold, so zero, and
the last vanishes by type. 
\end{proof}

The proof of Lemma~\ref{eform} leads to

\begin{prop}
  Assume $(N,J,h)$ is an almost Hermitian $6$-manifold with a half-integrable
  $\SU$-structure $(\omega,\psp)$. If the $\G$-structure \eqref{phi} on
  $N\times S^1$ has closed torsion $T$, then the K\"ahler form is an
  eigenvector of the Laplace-Beltrami operator on $N$:
  \begin{equation*}
    \Delta \omega=\tfrac{\lambda^2}{2}\omega.
  \end{equation*}
\end{prop}

\begin{proof}
  By~\eqref{T*T} $\Delta\omega = d^*d\omega =
  -\tfrac{\lambda^2}{4}\hodge_6\oo = \tfrac{\lambda^2}{2}\omega$, so
  $\omega$ is an eigenform with eigenvalue $\tfrac 12\lambda^2$.
\end{proof}

This fact partly motivates the reason for attracting the attention to the
following

\begin{defn}
  {\rm\cite{Cleyton-S:sG2}} A \emph{strong} $\G$-manifold with torsion
  (S$\G$T) is a $\G$-manifold with a closed skew-symmetric torsion form.
\end{defn}

Now, in the special case in which $N$ is a complex manifold equation
\eqref{dTd*T} seems very promising, for the derivative of the torsion becomes
\begin{equation*}
  dT=-d\hodge_6d\omega.
\end{equation*}
In other words if the $\G$-structure $\varphi$ is strong, then $N$ satisfies
\begin{equation*}
  dJd\omega=0,
\end{equation*}
recalling that the complex structure $J$ and the Hodge operator agree on
type $\{2,1\}$-forms.  The significance of this equation was uncovered
in~\cite{Alexandrov-I:vanishing}, then exploited in~\cite{Fino-PS:SKT},
together with the following fact of the utmost importance:

\begin{prop}{\rm\cite{Alexandrov-I:vanishing}} 
  A Hermitian manifold $(N,J,h)$ of dimension greater than four is SKT if
  and only if it is not balanced, i.e.~if the Lee form $\vartheta$ differs
  from zero.
\end{prop}
 
As $J\vartheta=d^*\omega$, the SKT notion is antithetical to that of
cosymplectic $\SU$-structures, and this conflict has dramatic consequences,
in the light of the fact (Fino et al.~\cite{Fino-PS:SKT}) that the holonomy
of an invariant Hermitian structure on a nilmanifold $N$ of dimension $n$
reduces to $\Lie{SU}(n)$ exactly when the metric is balanced. It basically
erases all chances for a strong $\G$-structure with torsion on $N\times
S^1$ to induce an SKT structure on $N$.

Recall that when $J$ is integrable, exterior differentiation is just
$\partial +\ol{\partial}$, so $d\hodge_6 d=2i\ol{\partial}\partial$ and the
SKT condition is the same as the $\partial\ol{\partial}$-lemma. But the
$6$-dimensional Lee form is essentially $\beta$, and asking $M^7$ to be
strong and $J$ a complex structure, forces $N$ to be K\"ahler
(cf.~Proposition~\ref{balanced}). Compare this to the result that a
compact, conformally balanced manifold with $\Lie{SU}(n)$ holonomy
satisfying $dJd\omega=0$ is in fact Calabi-Yau,
see~\cite{Papadopoulos:KT-HKT}.

This is confirmed by the fact that it is not possible to obtain strong $\G$
metrics from the Lie algebras we have considered, as no nilpotent Lie
algebra from Theorem~\ref{output} generates a structure with strong torsion
on products, unless $\lambda=0$. Because
\begin{equation*}
  d(JV_3)\subseteq\Lambda^2V_3\oplus (V_3\wedge JV_3)
\end{equation*}
in all instances, the common offending term in $dT$ has the expected
coefficient $-\lambda^2/4$ only in the degenerate case.

\providecommand{\bysame}{\leavevmode\hbox to3em{\hrulefill}\thinspace}
\providecommand{\MR}{\relax\ifhmode\unskip\space\fi MR }
\providecommand{\MRhref}[2]{%
  \href{http://www.ams.org/mathscinet-getitem?mr=#1}{#2}
}
\providecommand{\href}[2]{#2}

\end{document}